\documentclass{amsproc}

\usepackage{amssymb}

\usepackage{graphicx}

\usepackage[arrow,curve,matrix,tips,2cell]{xy}

\usepackage{hyperref}
\usepackage{pdfsync}
\usepackage{calc}
\usepackage{enumerate,amssymb}
\usepackage[arrow,curve,matrix,tips,2cell]{xy}

\DeclareMathAlphabet\EuR{U}{eur}{m}{n}
\SetMathAlphabet\EuR{bold}{U}{eur}{b}{n}

\makeindex  

\newtheorem{theorem}{Theorem}[section]

\theoremstyle{definition}
\newtheorem{definition}[theorem]{Definition}
\newtheorem{example}[theorem]{Example}

\newtheorem{condition}[theorem]{Condition}

\theoremstyle{remark}
\newtheorem{remark}[theorem]{Remark}

\numberwithin{equation}{section}

\newcommand{\calfin}{\mathcal{FIN}}

\newcommand{\IR}{{\mathbb R}}

\newcommand{\IZ}{{\mathbb Z}}

\newcommand{\curs}{\EuR}

\newcommand{\CHAINCOMPLEXES}{\curs{CHCOM}}

\newcommand{\FGINJ}{\curs{FGINJ}}
\newcommand{\MODULES}{\curs{MODULES}}
\newcommand{\Or}{\curs{Or}}

\newcommand{\coker}{\operatorname{coker}}

\newcommand{\sldrei}{Sl_{3}{\mathbb{Z}}}
\newcommand{\OO}{\mathcal{F}_{G}}
\begin{document}

\title{A  Survey of Computations of Bredon Cohomology}


\author{No\'e B\'arcenas }
\address{Centro de Ciencias Matem\'aticas. UNAM \\Ap.Postal 61-3 Xangari. Morelia, Michoac\'an MEXICO 58089}
\email{barcenas@matmor.unam.mx}
\thanks{ The author  thanks  the  organizers  of  the  session for  the  invitation  to participate at  the  AMS  joint  meeting  in  march  2022,    as  well as  support  from  DGAPA Project IN100221, the  DGAPA-UNAM Sabbatical Program and CONACYT through grant CB 217392.The  author  benefited  of  conversations  and collaboration  with  Mario  Vel\'asquez concerning   Bredon  cohomology  and  the  specific  examples  of  $\mathbb{S}l_{3}(\mathbb{Z})$. Special  Thanks  go  to Juanita  Claribel  Santiago, who made  figure 1. }

\subjclass[2010]{Primary:55N91. Secondary:55P91,55N25,20J05.}

\date{\today}

\begin{abstract}
We  present  an overview  of  computational  methods  for  Bredon cohomology with  a  special focus on infinite  groups. 
\end{abstract}

\maketitle



\section{Different Meanings  for Bredon Cohomology}

Bredon  cohomology  is  one  of the  most  prominent cohomology  theories  for  spaces  with  an  action of  a  group. 

We  reserve  the   notion  of equivariant ordinary cohomology -as understood  traditionally    and  in  this  volume- for  the cohomology  of  the  Borel  construction  of  a  space  with  an  action  of  a  compact Lie or  finite   group.  We  will  speak,  however,  of  Bredon  cohomology  as  an  equivariant  cohomology  theory  in  a  sense  to  be   defined  below.

Historically,  the  construction of  Bredon  cohomology  goes  back to   the  announcement \cite{bredonannouncement} and  the  extended  version  \cite{bredonlnm},  and it  is  strongly  based  on  the  notion  of  a $G$-CW  complex,  which  we  will  review in  this  note. 

It  is   this  relation  which   explains  the   use of  Bredon  cohomology   in the  study  of   finiteness  properties  in  group  cohomology \cite{lueckclassifying},  \cite{barcenasdegrijsepatchkoria}.  We  will  not  extend  in  the  discussion  of  this  subject, and  rather  refer  to  the  excellent  survey \cite{luecksurveyclassifying},   as well as to  \cite{nucinkis}, specifically  to the  Eilenberg-Ganea  problem  for  families which   is  phrased  in  terms  of  Bredon  cohomology of  classifying  spaces  for  families.

Simplicial versions  of Bredon cohomology were  provided  by Br\"ocker and  Illman \cite{broecker}, \cite{illman},  with the outcome   of  the  posibility  of  considering  actions  of (Hausdorff,  locally  compact) topological groups,  based  on   an  equivariant  version  of  simplicial  complexes, which mimics  the  definition  of $G$-CW  complexes.
 
Homotopical versions of  Bredon  cohomology,  which  even   allow  a  description  in the parametrized  equivariant  setting  are  described  in \cite{mukherjee} for  \emph{Na\"ive equivariant} cohomology  theories,  in the  sense  of  equivariant  homotopy  theory.  For  the  specific  case  of  Bredon  cohomology with local (twisted) coefficients in the  complex twisted  representation  group,  a  construction   has  been   written in  detail  in \cite{barcenasespinozauribevelasquez}, where  also  the  relation  to  the  \v{C}ech  versions   has  been  established  in  order  to  provide   a  description  of  the Segal/Atiyah-Hirzebruch spectral  sequence.

 \emph{Genuine  equivariant}   versions  of  Bredon  cohomology  have  been  considered  in  equivariant  homotopy  theory \cite{lewismaymclure}.  Classical  computations  of  these theories  include  seminal  work  of Gaunce Lewis \cite{gaunceprojective} (atributed to  unpublished  work  of Stong),  based  on  explicit  cellular structures  of  the  relevant  examples, the  $RO(G)$-graded  K\"unneth  and universal  coefficient  spectral  sequences \cite{lewismandell},  the  implicit  use  of  classifying  spaces  for  families,  such  as  the  Tits building  in  \cite{aronedwyerlesh},  the  explicit  use  of classifying  spaces for  families of  proper  subgroups \cite{kriznotes},
and  a  more  recent development  of a  variety of tools whose  interest goes back to  the  role  of   such  computations  in the  proof  of  the Kervaire  invariant  one  problem in \cite{hillhopkinsravenel}.   
 
 The  methods  include (without  the  intention  to  be  exhaustive in their  enumeration) those based  on  the  slice  filtration \cite{hillslice},   the  homotopy  fixed  point  spectral  sequences \cite{hollerkriz},  often in  combination  with   ad-hoc  cellular constructions \cite{krizlu},  as  well  as parametrized  homotopy  theory  considerations \cite{costenoblehudsontilson},  and  the  notion  of  freeness of \cite{hillfreeness}.  They  are  quoted  here  with  the  idea  of  giving  a representative  example  of  an  application  of each  kind  of  method.

 Finally,  a  genuine proper  equivariant   version  of Bredon  cohomology  has  been  defined in  terms  of  equivariant  homotopy theory  in example  3.2.16 in  \cite{degrijsehausmannlueckpatchkoriaschwede},  where  also  the  relation  of  the  extensions  of  the  gradings  from  $\mathbb{Z}$  to  the  equivariant  $KO^{0}$- theory  of  the  classifying  space  for  proper  actions (more  general  than  $RO(G)$-graded)  is  adressed.

 We  will  focus on  computational  methods  for  the  determination of  \emph{naive} versions  of  Bredon  cohomology,  with  an  emphasis  on  infinite  discrete  groups,  extending  the  content of  the  lecture delivered  at  the  AMS  sectional  meeting  with a  number  of  references  and  additional  examples  expanding  the exposition.

 \tableofcontents

 \section{Bredon Cohomology}

Let $G$ be a (possibly  infinite) discrete group.  A $G$-CW-complex is a CW-complex with a $G$-action permuting the cells and such that
if a cell is sent to itself, it is done by the identity map. We call the $G$-action proper if all cell stabilizers are finite
subgroups of $G$.

\begin{definition}
Recall  that  a  $G$-CW  complex structure  on  the  pair $(X,A)$  consists  of a  filtration of  the $G$-space $X=\cup_{-1\leq n } X_{n}$, $X_{-1}=\emptyset$,$X_{0}=A$  and  for   which  every  space  $X_{n}$ is inductively  obtained  from  the  previous  one   by  attaching  cells  in pushout  diagrams  of  the  form
$$\xymatrix{\coprod_{i} S^{n-1}\times G/H_{i} \ar[r] \ar[d] & X_{n-1} \ar[d] \\ \coprod_{i}D^{n}\times G/H_{i} \ar[r]& X_{n}}$$   
We  say  that a  proper  $G$-CW complex  is  finite  if  it  consists  of  a  finite  number  of  (orbits of) cells  $G/H\times  D^{n}$. 

\end{definition}

\begin{example}\label{examplecw}[Examples  of  $G$-CW  Decompositions]

\begin{itemize}
\item It  is  a  consequence  of  the  equivariant triangulation  theorem  that   there  exist equivariant  triangulations  of  smooth  manifolds  with  a  proper smooth  action  of  a Lie group  \cite{illmanproper},  \cite{illmantriangulation}.  Such  a  triangulation  produces a $G$-simplicial  complex.

 This notion  is  described  in  \cite{illmantriangulation}, section  5 as   a  triangulation  of  the  orbit  space  $X/G$ with  an  extra  compatibility condition with respect  to  the   quotient  map $\pi: X\to  X/G$,  namely:  the  inverse image  of  an $n$-dimensional  simplex $\Delta_{n}$ on $X/G$ is  a \emph{standard  equivariant  simplex}, denoted  $(\Delta_{n}, H_{0}, \ldots, H_{n})$,  which  is  a  quotient  of  the  product  of  a  free  $G$-orbit  of  the  standard $n$-dimensional  simplex  $\Delta_{n}\times G $,  where  pairs $(x,g)$ and  $(x,g^{'})$ are  identified if for  each  $k=0, \ldots, n$, the  element $x$  belongs  to  the boundary  of one  of  the  $k$-dimensional  simplices  in  $\Delta_{n}$, the   cosets for  the  compact  subgroup $H_{k}$  $gH_{k}$  and  $g^{'}H_{k}$  agree,  and  the  sequence of  compact  subgroups satisfy $H_{i+1}\subset H_{i}$. 
 
 Notice  that  the data  provided  by this identification   is  equivalent  to $G$-equivariant  maps  from  orbits of $i$-dimensional  simplices $G/H_{i}\times \Delta_{i}\to X$,  on  which   the  inclusions  of standard simplices $\Delta_{i}\subset \Delta_{i+1}$   
 are   made compatible with    inclusions  up to   $G$-conjugacy $H_{i+1} \subset H_{i}$. 

Such a  map  $\Delta_{i}\times G/H_{i}\to X$  can  be  identified  with  a cell. The  full  details  are  worked  out  in  \cite{illmantriangulation}, propositon  11.5, page 170.

 \item  A $G$-Absolute  Neighborhood Retract(G-ANR)  with   a proper  action  of  a  Lie  group,  or  more  generally,  a  locally  compact group has  the  homotopy  type  of  a  $G$-CW  complex. This is  is  a  consequence  of  the  slice  theorem \cite{palaisslice}.

 \item Equivariant cell  decomposition for  the  action of  $Sl_{2}(\mathbb{Z})$   on  the  hyperbolic  plane and  the   1- dimensional  deformation  retract.

The  group  $Sl_{2}(\mathbb{Z})$  acts  by  isometries  on  the  hyperbolic  plane. The  dual  of  the  Farey tesselation  provides  an  example  of  a  one- dimensional  $G$-CW  complex (the  Bass-Serre  tree for  $Sl_{2}(\mathbb{Z})$), consisting  of  two  orbits of  zero  dimensional cells of  type  $C_{4}$, $C_{6}$,  and  one  dimensional  cell of  type  $C_2$.  See figure 1. 

More   generally,  Li, L\"uck  and  Kasprovski constructed  in  \cite{kasprovskililueck} a  $G$-CW  structure  for  the  flag  complex  associated  to  a  group $G$  given as  a  graph  product
(Examples  include   right  angled  Artin  groups  and  right  angled  Coxeter  groups).  
  \begin{figure}\label{farey}\caption{Dual   of  the  Farey  triangulation.}
  
 \includegraphics[scale=0.5]{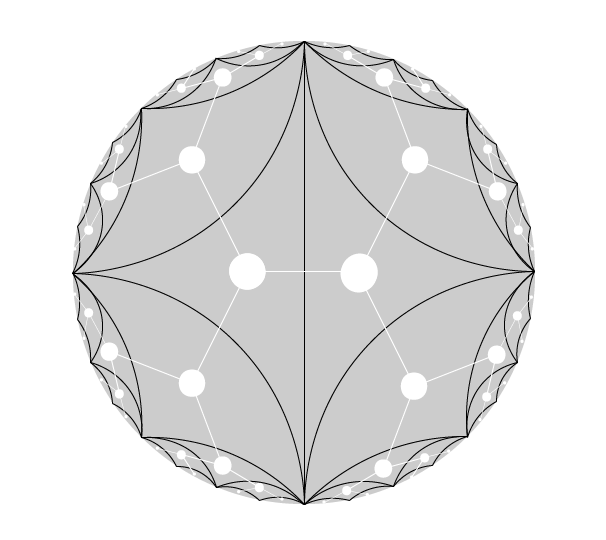}
\end{figure}  
 \item Equivariant  cell  decomposition  for the  action  of  $Sl_{3}(\mathbb{Z})$ on   the  homogeneous  space $Sl_{3}(\IR)/SO_{3}$. There  exists a  triangulation  of  the  quotient  of an  equivariant  deformation  retract of   this  homogeneous  space   described  in \cite{Soule},  but  also  in  
\cite{Sanchez-Garcia},  which  is  the  main  input  for   the  computations  of  twisted   equivariant $K$-Theory  and  $K$-homology in  \cite{barcenasvelasquez}, \cite{barcenasvelasquezrefrito}.

 Let  $Q$  be  the   space  of real, positive   definite $3\times  3$-square   matrices. Multiplication  
by  positive  scalars  gives  an  action  whose  quotient  space  $Q/ \mathbb{R}^{+}$  is   homotopy  equivalent  to   the  quotient $Sl_{3}(\IR) /SO_{3}/SL_{3}(\mathbb{Z})$.

We  describe its  orbit  space. Let  $C$  be  the  truncated  cube   of  $\mathbb{R}^{3}$  with  centre  $(0,0,0)$ and  side  length 
$2$,  truncated  at  the  vertices   $(1,1,-1), (1, -1,1),(-1,1,1)$  and  $(-1,-1,-1)$,  trough  the  mid-points  of  the   
corresponding sides. As  stated  in \cite{Soule}, every matrix  $A$  admits  a  representative of  the  form   

\[ \left( \begin{array}{ccc}
2 & z & y \\
z & 2 & x \\
y & x & 2 \end{array} \right)\] 

which  may  be  identified  with  the  corresponding  point  $(x,y,z)$  inside  the  truncated  cube. We  introduce  the  following  
notation  for  the  vertices of  the  cube: 
\[ \begin{array}{ccc}
 O=(0,0,0) & Q=(1,0,0)\\
M=(1,1,1) & N=(1,1,1/2)\\
M^{'}=(1,1,0) & N^{'}= (1,1/2, -1/2)\\
P=(2/3,2/3, -2/3) & \\
\end{array}
\]

Let $E$ be the  subspace  of  $C$ given  by the  points $(x,y,z)$ satisfying: 
$$\mid z\mid\geq  y \geq 1$$
$$ z-x-y+2\geq 2.$$ 

We  will  describe  an  equivariant  triangulation  of  the  space $E$. 

\begin{figure}[h]\label{figure1}
\caption{Triangulation  for  the  fundamental  region.}
\includegraphics[scale=15]{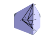}

\end{figure}

Notice that the elements of $\sldrei$  
$$q_1=\begin{pmatrix}1&0&0
                                    \\0&0&-1
                                    \\0&1&1\end{pmatrix}
    q_2=\begin{pmatrix}-1&0&0
                                    \\0&1&1
                                    \\0&0&-1\end{pmatrix}$$send  the  triangle  $(M,N,Q)$  to the  triangle  $(M^{'}, N^{'} Q)$ and  the  quadrilateral 
 
$(N,N^{'},M^{'},Q)$  to $(N^{'}, N, M^{'}, Q)$. Thus, the  following identification  must  be  performed in the  quotient: 
$M\cong M^{'}$, $N\cong N^{'}$, $QM\cong  QM^{'}$, $QN\cong QN^{'}$, $MN\cong  M^{'}N^{'}\cong M^{'}N$ and   
$QMN\cong QM^{'}N\cong QM^{'}N^{'}$.

Following \cite{Sanchez-Garcia} we now describe the orbits of cells and corresponding stabilizers. This can be found also in Theorem 2 of Soul\'e's article \cite{Soule}
(although we use a cellular structure instead of a simplicial one). We have changed
the chosen generators. We summarize the information  within 
Table \ref{table}. We use the following notations: $\{1\}$ denotes the trivial group, $C_n$ the cyclic group
of $n$ elements, $D_n$ the dihedral group with $2n$ elements and $S_n$ the Symmetric group of permutations on $n$ objects.

\vspace{0.5cm}

\begin{table}

\begin{tabular}{ccccccccc}

\hline
vertices & & & & &2-cells\\
\hline
$v_1$ & $O$ & $g_2$, $g_3$ & $S_4$& &$t_1$ & $OQM$ & $g_2$ & $C_2$\\
$v_2$ & $Q$ & $ g_4$, $g_5$& $D_6$&&$t_2$ & $QM'N$ & $g_1$ & $\{1\}$\\
$v_3$ & $M$ &  $g_6$, $g_7$&  $S_4$& &$t_3$ & $MN'P$ & $g_{12},g_{14}$ & $C_2\times C_2$\\
$v_4$ & $N$ & $ g_6$, $ g_8$ & $ D_4$&& $t_4$ & $OQN'P$ & $g_5$ & $C_2$\\
$v_5$ & $P$ & $g_5$ , $g_9$ & $ S_4$ && $t_5$ & $OMM'P$ & $g_6$ & $C_2$\\
\hline
edges&&&&& 3-cells\\
\hline
$e_1$ & $OQ$ & $g_2$, $g_5$ &  $C_2\times C_2$ && $T_1$& $g_1$ & $\{1\}$\\
$e_2$ & $OM$ & $g_6, g_{10}$ & $D_3$\\
$e_3$ & $OP$ & $ g_6, g_5$ & $D_3$\\
$e_4$ & $QM$ & $g_2$ & $C_2$\\
$e_5$ & $QN'$ & $g_5$ & $C_2$\\
$e_6$ & $MN$ & $g_6, g_{11}$ & $C_2\times C_2$\\
$e_7$ & $M'P$ & $g_6, g_{12}$ & $ D_4$\\
$e_8$& $ N'P$ & $ g_5, g_{13}$ & $ D_4$\\

\hline

\end{tabular}
\caption{Stabilizer  groups of  the  cells  under  the  $Sl_{3}(\mathbb{Z})$-action.\label{table}}
\end{table}

\vspace{0.5cm}

The first column is an enumeration of equivalence classes of cells; the second lists a representative of each class;
the third column gives generating elements for the stabilizer of the given representative; and the last one is the
isomorphism type of the stabilizer.  The
generating elements referred to above are

$$g_1=\begin{pmatrix}1&0&0\\
                                    0&1&0\\
                                    0&0&1\end{pmatrix}                    
     \hspace{1cm}g_2=\begin{pmatrix}-1&0&0\\
                                    0&0&-1\\
                                    0&-1&0\end{pmatrix}                               
      \hspace{1cm}g_3=\begin{pmatrix}0&0&1\\
                                    0&1&0\\
                                    -1&0&0\end{pmatrix}$$
                                    
$$g_4=\begin{pmatrix}-1&0&0\\
                                    0&1&1\\
                                    0&0&-1\end{pmatrix}                    
     \hspace{1cm}g_5=\begin{pmatrix}-1&0&0\\
                                    0&0&1\\
                                    0&1&0\end{pmatrix}                               
      \hspace{1cm}g_6=\begin{pmatrix}0&-1&0\\
                                    -1&0&0\\
                                    0&0&-1\end{pmatrix}$$
                                    
$$g_7=\begin{pmatrix}0&0&-1\\
                                    -1&0&0\\
                                    1&1&1\end{pmatrix}                    
     \hspace{1cm}g_8=\begin{pmatrix}-1&0&0\\
                                    0&1&0\\
                                    0&-1&-1\end{pmatrix}                               
      \hspace{1cm}g_9=\begin{pmatrix}0&0&-1\\
                                    -1&0&-1\\
                                    0&1&1\end{pmatrix}$$
                                    
$$g_{10}=\begin{pmatrix}0&0&-1\\
                                    0&-1&0\\
                                    -1&0&0\end{pmatrix}                    
     \hspace{1cm}g_{11}=\begin{pmatrix}-1&0&0\\
                                    0&-1&0\\
                                    1&1&1\end{pmatrix}                               
      \hspace{1cm}g_{12}=\begin{pmatrix}0&-1&-1\\
                                    0&-1&0\\
                                    -1&1&0\end{pmatrix}$$
                                    
$$g_{13}=\begin{pmatrix}0&1&1\\
                                    1&0&1\\
                                    0&0&-1\end{pmatrix}                    
     \hspace{1cm}g_{14}=\begin{pmatrix}-1&0&0\\
                                    -1&0&-1\\
                                    1&-1&0\end{pmatrix}$$                                    
                                    
We fix an orientation; namely, the ordering of the
vertices $O < Q < M < M' < N < N' < P$ induces an orientation in $E$ and also in   the  quotient  $\underline{B}\sldrei = E/\equiv$.  The cells coboundaries are  detemined  in  section  5  of \cite{barcenasvelasquez} and  include  restriction of  representations  and signs  coming  from  the prescribed  orientation  chosen  above. This  finishes  the  cell  structure  for  the  action.

  \end{itemize}

\end{example}

Bredon  cohomology  admits  a  description  in  terms  of functors  depending  on  the  orbit  category  of  a group. 

Similar  to  the  discussion  of  part  1   in  \ref{examplecw},  the  subconjugacy  relations  of subgroups,  say  from $H$ to  $K$   by  an  element  of  $G$, denoted  by  $g$ satisfying $gHg^{-1}=H^{'}\subset K $,  has  the  
consequence  that  there  exists  a $G$-equivariant  map $G/H\to G/K$ given  by  asigning  to  the  coset $g^{'}H$  the coset $gg^{'}g^{-1}K $ . 

 Such  equivariant  maps  determine  geometric  information  when  considered  as part  of  the intrinsically  given  data   in a  $G$-CW structure. In  particular,  such maps  origin  inclusions  of orbits  of  cells in  such  a  way  that  the  lowest  dimension  of  cells  is  associated  to  the  biggest  groups  in  terms  of  inclusion up  to  $G$-conjugacy.

 There  exist  two  important  categories  associated  to a  group,  which  are  relevant  to  the codification  of  such  relations  among the subgroups  of $G$,  the  orbit  category, and  the conjugation  homomorphism category. We  will  assume  for  the  rest  of  this  section  that  the  group $G$  is  discrete.

\begin{definition}[Orbit Category]
Denote  by  $\mathcal{O}_G$ the orbit category of $G$; a category with one object $G/H$ for each subgroup $H\subseteq G$ and where  morphisms are  given  by  $G$-equivariant maps.  There  exists  a  morphism $\phi:G/H\rightarrow G/K$ if and only if $H$ is conjugate in $G$  to a subgroup of $K$. More  generally,  given a  family  of   subgroups  $\mathcal{F}$,  which  is  closed  under  intersection and  conjugation,  we  can  form  the full subcategory $\mathcal{O}_{G}^{\mathcal{F}}$  where  the  objects  are  of  the  form $G/H$  with  $H$ in $\mathcal{F}$.  
\end{definition}

\begin{definition}[Conjugation Homomorphism Category]
Given a  group  $G$,  the conjugation-homomorphism  category $\mathcal{S}_{G}$  is  the  category  where  the  objects  are  subgroups  of $G$,  and   the  set  of  morphisms   between  two  objects $H$  and  $K$ is  the quotient, denoted  by  ${\rm Conhom}(H, K)/ {\rm Inn(K)} $  of  the  set  ${\rm conhom}(H,K)$ of  group  homomorphisms $\varphi:H\to K$ for which  there  exists  an element $g\in G$ such that  $\varphi$  is given  as  conjugation by $g$,  and  the  group  ${\rm  Inn}(K)$  of  inner  automorphisms of $K$ acts  by  composition.
 \end{definition}

There  exists  a  projection  functor ${\rm pr}:  \mathcal{O}_{G}\to 
\mathcal{S}_{G}$, which  asigns  to  each  orbit  $G/H$ the  subgroup $H$  and   to a  $G$- map $f: G/H\to G/K$ the   homomorphism $H\to K$   defined  as  conjugation  by  an element $g$  satisfying $gHg^{-1}\cong H^{'} \leq K$.

\begin{remark}\label{isomorphismcategories}

The  set ${\rm Mor}_ {\mathcal{S}_{G}}(H, K)$  is   isomorphic  to  the  quotient  of  the  action  of  the  centralizer of  the  subgroup  $H$  in $G$,   $C_{G}(H)=\{ g\in G \mid gh = hg\,   \text{for  all $h \in H$ } \}$  on  the   set  of $G$-equivariant  maps $\rm {Mor}_{\mathcal{O}_G}(G/H, G/K)$, where  an  element  of  the centralizer $g\in C_{G}(H)$  acts  by   composition  with the  right  multiplication $R_{g^{-1}}: G/H\to G/H$, $g^{'}H\mapsto g^{'}g^{-1}H$. 

The main  advantage  in considering  this  category  is  that  the   automorphism  group  of  a finite  group $H$ is  finite.

 In  the  orbit  category,  the automorphism  group  of  an object $G/H$ is  the   quotient  of  the  nomalizer  susbgroup  in  $G$ by  the  subgroup $H\cdot C_{H}(G)$  consisting  of  elements of  the  form  $hc$,  where  $h\in H$, and $c$ is  an  elment  of  the  centralizer  in  $G$ of  $H$.  We  will denote  this  group by $W_{G}(H)= N_{G}(H)/H\cdot C_{H}(G)$.

 Both the  orbit  category  and the conjugation-homomorphism category  are EI- categories,  in  the  sense  that  every  endomorphism  of  an  object  is  invertible.

\end{remark}

\begin{example} \label{exampleorbits}  [Orbit  Categories  for  Infinite Groups]
\begin{itemize}
\item The  orbit  category   for  the  family  of finite subgroups  for  the  group $Sl_{2}(\mathbb{Z})$  has    three  objects: $Sl_{2}(\mathbb{Z})/C_{6}$, $Sl_{2}(\mathbb{Z})/C_{4}$,   and $Sl_{2}(\mathbb{Z})/C_{2}$.  There  exist  two  $G$- equivariant  maps $ Sl_{2}(\mathbb{Z})/C_{2}\to Sl_{2}(\mathbb{Z})/C_{6}$ and $Sl_{2}(\mathbb{Z})/C_{2}\to Sl_{2}(\mathbb{Z})/C_{4}$.

\item The  triangulation  of  the  quotient  of  the action  of $Sl_{3}(\mathbb{Z})$  on  the  space constructed  by Soul\'e  and  discused  in \ref{examplecw}, part  4,  has   the  consequence  of a complete  description  of  the orbit  category  for   finite  stabilizer  subgroups in   the  group: There  exist eight  maximal  finite  subgroups, which  are  the  stabilizers  of  the  zero  dimensional  cells,  there   are  as  many  morphisms  between  them  as  the  one  dimensional  cells of the  triangulation  having  as  edges  the  vertices, the  composition  of  pairs  of  such  morphisms  are  related  by   the  obvious  rule  given by the   two dimensional  cells,  and  finally,  there  exists  a  unique  composition  of length three  in the  orbit  category  for  finite  isotropy  subgroups of $Sl_{3}(\mathbb{Z})$.  See  the  table  at  the  end  of  the  previous section.

\end{itemize}

\end{example}

\begin{definition}[The $\mathcal{O}_{G}$-Chain Complex Associated  to a  $G$-CW complex]
Let  $X$  be  a  $G$-CW complex. Denote  by $\CHAINCOMPLEXES$ the category of chain complexes  over a ring $R$, and  morphisms given by  chain  maps.
 
 The  contravariant   functor  $\underline{C}_{*}(X):\mathcal{O}_G\to  \mathbb{Z}-\CHAINCOMPLEXES  $  assigns  to every  object  $G/H$     the  cellular $\mathbb{Z}$- chain  complex   of  the $H$-fixed point  complex    $ \underline{C}_{*}(X^ {H})\cong C_{*}({\rm  Map  }_{G}(G/H, X))$  with  respect  to  the  cellular  boundary  maps $\underline{\partial}_{*} $. 
\end{definition}

We  will  use  homological  algebra  to  define Bredon  homology  and  cohomology  functors. 

A  contravariant Bredon Module  is  a  contravariant  functor $ N:\mathcal{O}_{G}^{\mathcal{F}}\to \mathbb{Z}-\MODULES$,  where  $\OO$  is   the  full  subcategory  of the  orbit  category  of  $G$,    $\mathcal{O}_G$  generated   by  the  objects $G/H$ for  a  family of  subgroups  $H\in \OO$.

Given  a  contravariant Bredon  module $M$, the  Bredon    cochain   complex $C_G^*(X;M)$ is  defined  as the   abelian  group   of  natural  transformations   of  functors  defined  on  the  orbit  category $\underline{C}_{*}(X) \to  M$. In  symbols, 

$$C_G^n(X;M)=Hom_{\mathcal{F}_G}(\underline{C}_n(X),M)$$
Where  $\mathcal{F}_G$ is  a  subcategory  containing  the  isotropy  groups  of  $X$.

Given a  set  $\{e_{\lambda}\}$ of   orbit   representatives of  the n-cells of  the  $G$-CW  complex  $X$,  and isotropy  subgroups  $S_{\lambda}$  of  the  cells  $e_{\lambda}$,  it  is  a direct  consequence of  the  Yoneda  lemma that  the  abelian  groups $C_G^n(X,M)$  satisfy:
 
 $$C_G^n(X,M)= \underset{\lambda}{\bigoplus }Hom_{\mathbb{Z}}(\mathbb{Z}[e_{\lambda}], M(G/S_{\lambda}))$$   
 with  one  summand  for  each  orbit representative  $e_\lambda$.
They  afford  a differential $\delta^n:C_G^n(X,M)\to C_G^{n-1}(X,M)$ determined  by  $\underline{\partial}_*$ and  maps $M(\phi):  M(G/S_\mu)\to M( G/ S_\lambda )$  for  morphisms  $\phi:G/S_\lambda \to  G/S_\mu$.  

 For   a  covariant functor  $$N:\OO\to  \mathbb{Z}-\MODULES,$$  the chain  complex 

$$ C_{*}^ {\mathcal{O}_G}= \underline{C}^{n}(X) \underset{\OO}{\otimes}N =\underset{\lambda}{\bigoplus}\mathbb{Z}[e_{\lambda}]\otimes  N(G/S_{\lambda})$$ 

admits  differentials $\delta_* =\partial_* \otimes N(\phi)$  for    morphisms   $\phi:  G/ S_{\lambda}\to  G/S_{\mu}$ in $\OO$.

\begin{definition}[Bredon Homology]
  The  Bredon  homology   groups  with  coefficients  in $N$ ,  denoted  by  $H_{*}^{\mathcal{O}_G}(X, N)$,   are  defined  as  the  homology  groups   of  the  chain  complex $\big( C_{*}^{\mathcal{O}_{G}}(X, N), \delta_* \big)$
\end{definition}

\begin{definition}[Bredon cohomology] 
The  Bredon  cohomology  groups   with  coefficients  in  $M$, denoted  by  $H^{*}_{\mathcal{O}_G} (X,  M)$    are  the  cohomology   groups  of  the  cochain  complex  $\big (C_{G}^*(X, M), \delta^ * \big )$. 

\end{definition}
    
\begin{remark}[Bredon Cohomology  in  Terms  of  the  Conjugation Homomorphism Category]
Let  $X$  be  a  $G$-CW  complex with finite  cell  stabilizers. 
Let  $R$ be  a  ring,  and let $M$ be   a contravariant  Bredon  module  taking  values  on  $R$-Modules. 

We  will consider    the  contravariant  functor defined  on  the conjugation  homomorphism  category $\mathcal{S}_{G}$  and  taking  values  in the  category  of chain  complexes,  where  we  asign  to  a  subgroup $H$ the  cellular  chain complex   of  the  quotient  space $X^{H}/ C_{G}(H)$. 
We  will  denote by  $C_{*}^{\mathcal{S}_{G}}$ the obtained functor,  and  define  for  every  contravariant  functor  defined  on $\mathcal{S}_{G}$ and  values  on  the  category of  $R$-modules  the  cochain  complex  of  natural  transformations  from   $C_{*}^{\mathcal{S}_{G}}$  to  $M$  as 
$Hom_{\mathcal{S}_{G}}(C_{*}^{\mathcal{S}_{G}}, M)$. We  can  apply  the cohomology  funtor  with  respect  to the  cellullar  cochain  maps  and  we  will  denote  the   obtained  modules   by $ H^{n}_{\mathcal{S}_{G}}( C_{*}^{\mathcal{S}_{G}}, M)$.

As  a  consequence of  remark \ref{isomorphismcategories},  the  projection  functor ${\rm pr}: \mathcal{O}_{G}\to \mathcal{S}_{G}$ induces  for  every  contravariant  functor $M$ defined  on the   conjugation  isomorphism  category and  taking  values  in  $R$-modules  a  bijective  pair    of   natural  transformations 
$$pr^{*}: {\rm pr}^{*} M \to M,	 $$ 
 between  the  composition  of  functors, and 
$$pr_{*}: C_{*}^{\mathcal{O}_{G}}\to C_{*}^{\mathcal{S}_{G}}  $$
between  the  chain complexes,  which  produce  an  isomorphism  on  the  level  of  cohomology  groups 

$$H^{n}_{\mathcal{O}_{G}}(X, {\rm pr}^{*}M) \cong H^{n}_{\mathcal{S}_{G}}( C_{*}^{\mathcal{S}_{G}}, M)$$

\end{remark}

\begin{example}[Examples  of  Bredon  Modules]
We  will  now  enumerate  some important  examples  of  Bredon  modules. 

\begin{enumerate}
\item Constant  coefficients.  Given a  fixed  $R$-Module  $M$,   we  consider the  constant  functor   with value  $M$  for  each  orbit. 
\item  Free Bredon  functors.  A contravariant  functor  defined  on  the  orbit  category  is said  to  be    free  if  it  is  a  direct  sum  of  representable  contravariant functors;  that  is,  there  exist  a  number  of   objects  $G/H_{\alpha}$ such  that  the  functor  is given  as 
$$ \bigoplus_{G/H_{\alpha}} R[{\rm Mor}(\quad , G/H_{\alpha})]. $$

Chain  complexes   associated  to  $G$-CW complexes provide    examples  of free  contravariant  Bredon modules.

\item The  complex  representation  ring as  a Bredon  module. 

We  restrict  to  the family  of  finite  subgroups  of  the   possibly  infinite  group  $G$,  and  associate   to  the  object in the orbit  category $G/H$  with  a  finite   subgroup $H$,  the  complex  representation  ring $R_{\mathbb{C}}(H)$.  Recall  that the  subjacent abelian  group is  free in  the  set  of  conjugacy  classes  of  $H$,  but  it  is  not  free  as  a  functor  over  the orbit  category  as   defined  in  the  example  above. 

On  the complex  representation  ring  we  can consider  a  covariant  structure,  assigning  to a $G$-map $G/H\to G/K$  the  induction  homomorphism $R_{\mathbb{C}}(H)\to R_{\mathbb{C}}(K)$  constructed  as  follows. Let  $H^{'}\leq K$ be  the   subgroup  of $K$ which  is   conjugated  to $H$. 

Let  now $V$ be   a  representation of  $H$, and  consider   it  as  an  $H^{'}$- representation. Consider  now  the  $K$-representation   
$V\otimes_{\mathbb{C}}\mathbb{C}(K)$. 

There  exists a  contravariant  structure  on  the  representation  ring,  defined  by  assigning  to  a  morphism $G/H\to G/K$  the   homomorphism $R_{\mathbb{C}}(K)\to R_{\mathbb{C}}(H^{'})\cong R_{\mathbb{C}}(H)$ given  by  restriction.

\item Twisted  complex  representation  rings

Let  $H$  be a   finite  group and  $V$  be  a  complex  vector  space. Given  a   cocycle     $\alpha:H\times H\to S^{1}$  representing  a  class  in  $ H^{2}(H,S^{1})\cong H^{3}(H,\mathbb{Z})$, an $\alpha$-twisted  representation  is  a   function $P:H\to Gl(V)$  satisfying: 
$$P(e)={\rm 1}$$  
$$P(x)P(y)=\alpha(x,y)P(xy)$$

The  isomorphism  type  of an $\alpha$-twisted  representation  only  depends  on the  cohomology  class  in  $H^{2}(H,S^{1})$.  

The  Bredon  module  structure   is  similar  as  in  the  previous example.

\begin{definition}
Let  $H$  be  a  finite  group  and   $\alpha:H\times H\to S^{1}$  be  a cocycle  representing a class  in  $ H^{2}(H,S^{1})\cong H^{3}(H,\mathbb{Z})$.  The  $\alpha$-twisted  representation   group of $H$, denoted  by  $^{\alpha}\mathcal{R}(H)$  is   the  Grothendieck  group  of  isomorphism classes  of complex, $\alpha$-twisted representations with direct  sum as binary  operation. 
\end{definition}

Let   $H$  be  a  finite  group.  Given  a   cocycle  $ \alpha \in H^{2}(H,S^{1})$ representing  a  torsion class of  order  $n$,  the   normalization  procedure gives  a  cocycle $\beta$  cohomologous  to  $\alpha$  such that  $\beta: H\times  H\to  S^ {1}$ takes  values  in  the  subgroup  $\mathbb{Z}/n\subset S^{1}$   generated  by  a  primitive  $n$-th  rooth  of  unity $\eta$.  Associated  to  a  normalized  cocycle, there  exists   a   central  extension  
$$1\to  \mathbb{Z}/n\to H^*\to H\to  1$$
with  the  property  that  any   twisted   representation  of $H$  is  a linear  representation  of  $H^*$, with   the  additional  property  that  $\mathbb{Z}/n$  acts  by   multiplication  with  $\eta$. Such a  group  is called  a  Schur covering  group for  $H$.

\item Burnside  Ring.
We  restrict  again  to  the  family  of  finite  groups for  the  isotropy  of  objects  in  the  orbit  category.  Given  a  finite  group $H$, the  Burnside  ring $A(H)$ is  defined  to  be  the  Grothendieck  ring  of  the  set  of  isomorphism classes  of  finite  sets  with an  action  of  $H$. 

Similar  to  the  representation  ring, there  exist  two   structures  on  the  Burnside  Ring:  one   covariant,  given  by  induction  of  actions, and  a contravariant  one,  which  is  defined  by  restriction.

\end{enumerate}
\end{example}

\begin{example}[ Computation  of  Bredon Homology  of  $Sl_{2}(\mathbb{Z})$ from  the  Definition.] 
Let  us  consider the  complex  representation ring  as  Bredon  module.  We  will  restrict  to the  family  of  finite  subgroups here. 
Directly  from  the 1- dimensional  cellular  structure  for  the   (Bass Serre)  graph  given  as  de  dual  of  the  Farey triangulation, we  obtain  a  Chain  complex  computing  the  Bredon  homology  where  the  depicted  map  is  $d_{1}$. 

$$0 \to \mathcal{R}(C_{2}) \overset{\big ({\rm ind}_{C_{2}}^{C_{6}}),    -{\rm ind}_{C_{2}}^{C_{4}} \big )}{\to}\mathcal{R}(C_{6})\oplus \mathcal{R}(C_{4})\to 0 .$$

Notice  that  the  rank  of  the  abelian  group  on the  left  is  2, the  rank  of  the group  on  the  right  is  6+4=10,  and  the map is  injective. 

Thus,  if  we denote  by  $T$ the  $1$-dimensional  $Sl_{2}(\mathbb{Z})$-  complex  obtained.

$$ {H_{0}}^{\mathcal{O}_{Sl_{2}(\mathbb{Z})}}( T, \mathcal{R})=\mathbb{Z}^{8}, $$
and  all  other  Bredon  homology  groups  are  zero. 

\end{example}

\begin{remark}[Further Accesible Examples of  low  Dimension]
We  mention  for  the  sake of  utility    the   following  illustrative  examples  of  computations  of  Bredon  cohomology  in low  dimension: \cite{rumi}, \cite{flores}. 
\end{remark}

\begin{remark}[Computations for $SL_{3}(\mathbb{Z})$ Based  on the $G$-CW  Decomposition.]

The  triangulation  proposed  by  Soul\'e   has  been  extensively  exploited  for  computations  of Bredon cohomology. 

The   first  computation of  the  Bredon  chain and  cochain  complex  was   done in \cite{Sanchez-Garcia}. In  detail, the  determination  of stabilizers, and  conjugacy  relations, as  well  as  the  differentials  have  been  used  as input  for  computations   of  Bredon  cohomology    with  coefficients  in  complex  representations  in \cite{Sanchez-Garcia},  for   coefficients   in  twisted  representations  in \cite{barcenasvelasquez}, \cite{barcenasvelasquezrefrito},  and  more  recently,  with  the  coefficients  of  equivariant  real  $K$- Theory  in \cite{hughes}. 

The  output  is  that  the  spectral  sequence  of  Atiyah-Hirzebruch type  collapses, and the  computation  amounts  to  a  computation of  the  left  hand  side  of  the  Baum-Connes  conjecture.   

\end{remark}

\begin{remark}[Equivariant Obstruction Theory on $G$-CW Complexes and  Bredon Cohomology]
One  of  the  first  succeses  of Bredon cohomology  was  the   development of  equivariant  obstrution theory. Consider  the  equivariant  version  of  the obstruction problem: Given $G$-CW  complexes  $X$ and  $Y$,  and a   $G$-map $f: X_{n-1}\to Y$  defined  on  the $n-1$-skeleton, the  main  theorem  states  that  $f$  can  be redefined  to be $G$-homotopic  over  the $n-2$  skeleton  to a  map  which  extends over  an  $n+1$-dimensional additional orbit   of a cell $G/H\times D^{n+1}$  
 if  an  obstruction class vanishes 
$$\mathfrak{o}(f)\in H^{n}_{\mathcal{O}_G}(X, \pi_{n+1}(Y^{H})).$$
See \cite{bredonlnm} Chapter II and \cite{blagojevic} for  the  use of  equivariant  obstruction theory in  combinatorial  geometry in the specific  case of  free  actions. 
 
\end{remark}

   \begin{remark}[Miscellaneous ]
   We  mention the  following  recent  appearances  of Bredon  Cohoology in adjacent areas. 
\begin{itemize}
\item Lower  bounds  for  topological  complexity  of Eilenberg-Maclane spaces  have  been  obtained  with  computations  of  Bredon  cohomology  in \cite{farber}. 
\item There  exists a  framework  for  the  study  of  Hecke  operators  acting  on  Bredon  cohomology to  obtain a  Hecke  action on  the   reduced  $C^*$ algebra  of a Bianchi group  in \cite{velasquezhecke}. 

\item  There  exists  a  computation of a  motivic  version  of  Bredon  Cohomology  in  \cite{heller}. 
\end{itemize}

 \end{remark}

\begin{remark}[Algebraic  properties  of  the  Abelian  Category  of  Bredon  Modules  and   its Objects.]

We briefly  mentioned  in the  previous  example  that (by  definition of  free  object)  the  chain  complex  associated  to  a $G$-CW  complex  provides   an  example of  a  free  functor  over  the  orbit  category. We  will  examinate  some  algebraic  properties  of  the  categories  of modules  and chain  categories. 

The  category  where  the  objects  are  functors  on  the   orbit  category  in  a  category of  modules  is  an  abelian  category,  where  a morphism  is  a  natural  transformation  between  them.
A pair  of  consecutive  morphisms  is  said  to  be  exact  in the  middle if   it  is  exact  on  every  object. 

The notions  of  projective  and  injective  module can  be  given as  usual in terms  of  Hom  functors,  and  free and  projective  resolutions  exist  as  a  consequence  of  the  Yoneda Lemma. 

See \cite{nucinkis} for  an  introduction  to  Bredon homological  algebra in  connection  with group  cohomology, and  specifically  finiteness properties. 
 
\end{remark}

\begin{remark}[Computation of  Bredon  cohomology  in  practice]

The  coefficient  systems  considered  in  this  note  yield chain  complexes,  respectively  cochain complexes    of   free  abelian  groups with  preferred  bases  to  compute  both  Bredon  homology  and  cohomology. 

Notice that if we have a complex of free abelian groups
$$\cdots\rightarrow\IZ^{\oplus n}\xrightarrow{f}\IZ^{\oplus m}\xrightarrow{g}\IZ^{\oplus k}\rightarrow\cdots$$
with $f$ and $g$ represented by matrices $A$ and $B$ for some fixed basis, then the homology at $\IZ^{\oplus m}$ is
$$ker(g)/im(f)\cong\IZ/d_1\IZ\oplus\cdots\oplus\IZ/d_s\IZ\oplus\IZ^{\oplus(m-s-r)},$$
where $r = rank (B)$ and $d_1,\ldots, d_s$ are the elementary divisors of $A$.

This  has  a consecuence the \emph{Torsion Freeness  Criterion}  for  Bredon  homology,  see  Theorem 5.2 in  page 1496 of \cite{lafont}: if there exists a base transformation such that all minors in all vertex blocks are in the set $\{- 1, 0, 1 \}$,  then  the  zeroth graded Bredon  cohomology  is torsion-free. Effective  approaches  for  implementing  Bredon  cohomology  in  computations (mainly  in  GAP) include  \cite{sanchezgarciathesis},  \cite{rahm}. 

\end{remark}

\subsection{Comparison to  Other  Constructions}

While  the simplicial  and  cellular  versions  of  Bredon cohomology  are  compared  by Illman in \cite{illman}, the  \v{C}ech versions  is   seen  to  agree for  a  proper  $G$-ANR with the \emph{naive}-homotopical  versions  of  Bredon  cohomology  in the  more  general  parametrized  setting in Theorem  5.3  of \cite{barcenasespinozauribevelasquez}. The  theorem  is  phrased    for  the  specific  example  of  (twisted) complex  representation  ring,  but the  argument  holds  for  any Bredon  module.

\begin{remark}[Contrast  to  Genuine  Equivariant  Homotopical  Versions  of  Bredon  Cohomology]

There  exist versions  of  Bredon cohomology  obtained  by  considering  a  \emph{genuine equivariant} version  of  the  Eilenberg  Maclane  spectrum  in  a category  of  equivariant  spectra. 

 Some  examples  of  these  computations  include 
\cite{lewismandell}, \cite{aronedwyerlesh} \cite{kriznotes}
\cite{hillhopkinsravenel}, \cite{hillslice},  \cite{hollerkriz} \cite{krizlu}, \cite{costenoblehudsontilson}. 
 
We  refer  to   example 3.2.16 in  \cite{degrijsehausmannlueckpatchkoriaschwede} for  the  relation  of naive  graded  equivariant  cohomology  theories to   extensions  of  them  to genuine  ones,  which  apply  in  particular  to Bredon  cohomology.  
\end{remark}

\subsection{Computations  Based  on  the  Algebraic  Properties  of  the  Category  of Mackey  Functors}

We  will  present  now  a  decomposition of  Bredon  cohomology  in   terms  of  information  concerning automorphism  groups  of  each object for  the  Bredon module. The  relation  is   crucial  to   decompositions   which  refine  equivariant  cohomological Chern characters.

References  for  this  section  are  \cite{slominska}, \cite{lueckequivariantcohmological}, \cite{lueckcreelle}, from  where  the  totality  of  arguments and  definitions are  extracted with  no claim of  originality.

It turns  out  that there  exists  a  close  relation between  the  algebraic  properties  of  being  injective in  the  category  of  functors  over  the orbit  category, and  the possibility  of  the decomposition  in terms  of the  stabiizer  groups  associated  to  the  isomorphism  classes  of  objects. 

Consider  for  this,  in  order  to  fix  notation,  the  inclusion  of  an  object  in  the  orbit  category  $i_{G/H}: G/H\to  \mathcal{O}_G$, and  recall  that  the  automorphism  group of $G/H$  is   $W_{H}(G)$. 

Let  $R$ be  a  ring, and  let $M$  be a  contravariant  Bredon  module. 

There  is  the \emph {restriction} functor $i^{*}$  which  restricts $M$  to  the  $R$-module  $i^{*}(M)=M(G/H)$,  and  which gives  an $R[W_{G}(H)]$-module  structure  by  considering  the full subcateory  of  the  orbit  category  with  exactly  one object, $G/H$,  and  the  automorphisms  of  $G/H$ as morphisms.

Given  an $R[W_{H}(G)]$- module $M$, the   \emph{induced module} $i_{*}M$, is  the $\mathcal{O}_G$- R  module defined  as quotient of  the  tensor  with  the  contravariant free  object
$$G/K \mapsto M(G/H)\otimes R[{\rm Mor}(G/K, G/H)], $$
where  we declare  equivalent  to  zero  the  submodule of  $M(G/H)\otimes R[{\rm Mor}(G/K, G/H)]$  generated  by  elements  of  the  form $m\otimes f_{*}(x)- f^{*}m\otimes x$,  for  all  $f\in W_{G}(H)$.

Given  a $R[W_{G}(H)]$- module $N$,  the \emph{coinduction}   functor assigns (contravariantly!) to an object $G/K$ in the  orbit  category the $R$-module 
$${i_{G/H}}_{!}= {\rm Hom}_{R[W_{H}(G)]} (R[{\rm mor}_{\mathcal{O}_{G}}(G/H, G/K)] , N) ).$$

Notice  the  adjunctions for every pair consisting  of  an $\mathcal{O}_{G}$- functor  $M$ and  $N$  a  $R[W_{G}(H)]$-module. 

$$ {\rm hom}_{R[W_{G}(H)]} ( {i_{G/H}}^{*}  M, N)\cong  {\rm hom}_{\mathcal{O}_{G}}(M, {i_{G/H}}_{!} N)  $$

$$ {\rm hom}_{\mathcal{O}_{G}}(  {i_{G/H}}_{*} N,  M) \cong {\rm hom}_{R[W_{H}(G)]}(N, i^{*}_{G/H}(M)).$$

\begin{definition}[Projective  and  Injective Splitting  Functors]
 Let $G/H$ be  an object  in  the  Orbit  category. 
\begin{enumerate}

\item 
 The  projective splitting  functor  $S_{G/H}$ associates  to   a  contravariant  functor defined  over  the  orbit  category  $M$, the $ R[W_{G}(H)]$  module defined  as the cokernel  of  the  map 
 $$\underset{\underset{\text{$f$ not an  isomorphism}}{f: G/H\to G/K}}{\prod}M(G/K)\to M(G/H). $$

\item  The  injective  splitting  functor $T_{G/H}$ associates  to  a  covariant  functor  defined  over  the  orbit category $N$,  the $R[W_{G}(H)]$-module  defined as  the  kernel  of  the  map 
$$\underset{\underset{\text{$f$ not an  isomorphism}}{f: G/H\to G/K}}{\bigoplus} N(G/H)\to N(G/K).$$

\end{enumerate}
\end{definition}

The  projective splitting  functor  comes  with  a  canonical  projection  $M(G/H)\to S_{G/H}(M)$. Given  any $R[W_{G}(H)]$-section,  the  inclusion  of  the  object $G/H$  into  the orbit  category  produces  a natural transformation ${i_{(G/H)}}_{*}S_{G/H}(M) \to{ i_{G/H}}_{*}(M)$. 
Dually,  the  projective  splitting  functor   comes  with  a  canonical  injection 
$T_{G/H} M \to M(G/H)$, an  any  $W_{G}(H)$- retraction $M(G/H)\to T_{G/H}M$ produces a natural  transformation  
$$ {i_{G/H}}_{!}M(G/H) \to  {i_{G/H}}_{*}T_{G/H}M.$$ 

For  an  object  $G/H$ in either  one  of  the  categories  $W_{G}(H)$ or $\mathcal{O}_{G}$, the  length  of  an  object $G/H$ is  the  supremum  of  all $l$  for  which  there  exists a  sequence  of  morphisms 
$G/H_{0}\to  G/H_{1}\to \ldots \to G/H_{l}$ with  $G/H_{l}=G/H$  and   none  of the  morphism  is an isomorphism,  dually, the  colength  is  the  supremum over all  $l$  for  which  there  exists  a  sequence  $G/H_{0}\to  G/H_{1}\to \ldots \to G/H_{l}$ with   $G/H_{0}= G/H$, and    none  of  the  morphisms  is  an  isomorphism.
A  category  has  finite  length,  respectively  finite  colength,  if  each  object  has  finite  length  or  colength.

The  Structure  Theorem \cite{lueckequivariantcohmological}, 2.2  in  page  1035 reads  as  follows: 
\begin{theorem} \label{structuretheorem}
\begin{enumerate}
\item Suppose  that  $\mathcal{O}_{G}$ has  finite  colength,  and  that  $M$  is  a  covariant Bredon  functor with  the  property  that $S_{G/H}M$  is a projective $RW_{H}(G)]$- module for  each object $G/H$.   Let  $\sigma_{G/H}: S_{G/H}M\to M$ be  an $R[W_{H}(G)]$- section  for  the canonical  projection,  and  consider  the  map  of Bredon  functors 

\begin{multline*}
\mu(M): \\ \underset{G/H\in {\rm Iso}(\mathcal{O}_G)}{\bigoplus}{i_{G/H}}_{*} S_{G/H} M \overset{ \underset{G/H\in {\rm Iso}(\mathcal{O}_G)}{\bigoplus {i_{G/H}}_{*} \sigma_{G/H}M}   }{\longrightarrow} \\  {i_{G/H}}_{*} M(G/H) \overset{\underset{G/H\in {\rm Iso}(\mathcal{O}_G)} {\bigoplus \alpha_{G/H}}} {\longrightarrow}M.
\end{multline*}
Where $\alpha: {i_{G/H}}_{*}M(G/H)= {i_{G/H}}_{*} i_{G/H}^{*}M \to M$  is  the  adjoint  of  the  identity. 

The  map  is  always surjective.  It  is  bijective  if  and  only  if  $M$ is  a projective Bredon module. 

\item Suppose  that  $\mathcal{O}_{G}$ has  finite length. Let  $M$ be  a  contravariant $R_{\mathbb{O}_G}$-module  such  that  the $R[W_{G}(H)]$- module $M(G/H)$ is  injective  for  every $G/H$. Let $\rho_{G/K}:M(G/K)\to T_{G/K} M$ be  an $R[W_{G}(H)]$- retraction  of  the canonical  injection $T_{G/H}M \to M$ and  consider  the   natural  transformation

\begin{multline*}
 \nu(M): \\ M \overset{\prod_{G/K \in {\rm Iso(\mathcal{O}_{G})}}\beta_{G/K}}{\longrightarrow} \\ \underset{G/K \in {\rm Iso}(\mathcal{O}_{G})} {\prod}	{i_{G/H}}_{!}M(G/K)
\overset{\underset{{G/K \in {\rm Iso}(\mathcal{O}_{G})}}{\prod} {i_{G/K}}_{!}\rho_{G/K}} {\longrightarrow} {i_{G/K}}_{*}T_{G/K}M.  
\end{multline*}
The  map  is  always  injective. It  is  bijective  if $M$ is  an injective Bredon  module.

\end{enumerate}

\end{theorem} 

In the  important  example  of   complex representation  rings, we  noticed  the   fact  that  there  exists a  covariant  and contravariant structure on  fuctors  which  agree  on  objects. 

On the  other  hand  side,  the  second  part  of Theorem   \ref{structuretheorem},  a  condition  appears  as  for  the  map $\nu_{G/H}$  to  be  surjective,  which  is  equivalent  to  the  fact  that  the Bredon  module  is  injective.  This  has  the   consequence  for  a  contravariant  functor $M$ that   the   composition  with   the  projection functor  to  the  conjugation  homomomorphism  category $pr^*(M):\mathcal{S_{G}}\to {\rm R-Mod}$  gives    an  injective  $R$-module after evaluation  on  each  object. 

The  following  notion  is  an  equivalent   characterization of this  property,  which   and  will  be  the   most  relevant algebraic  tool  to  the  rational  computation  of  Bredon  cohomology. 

 \begin{definition}
 Let  $ \FGINJ$ be  the  category  of  finitely  generated groups  and  injective  group  homomorphisms.  Let $M^*, M_*$ be  a  bifunctor to  the  category  of $R$-modules;  that is,  a  pair  consisting  of  a  contravariant  functor $M^*$  and  a  covariant functor  $M_*$ agreeeing  on  objects.  We  will  denote by   ${\rm  ind}f$  the  covariantly  induced   homomorphism,  and   by ${\rm res} f$    the contravariantly  induced  homomorphism.  For  inclusions  of  a  subgroup $H\to G$,  we  will  write ${\rm res}^{H}_{G}$  and  ${\rm ind}_{H}^{G}$. 
 
  M  is  said  to  be  a  Mackey functor  if 
 \begin{itemize}
 \item For  an  inner  automorphism $c(g):G\to G$,  we  have $M_{*}(c_g): M(G)\to M(G)  $ is  the  identity.
 \item For   an  isomorphism  of  groups $f:G\overset{\cong}{\to} H$,  the  composites  ${\rm res} f\circ {\rm ind} f$ and ${\rm ind} f\circ {\rm res}f$ are  the  identity. 
 \item Double  coset  formula.  For  two  subgroups $H, K\subset G$, 
 $$ {\rm res}^{K}_{G}\circ {\rm ind }_{H}^{G}=\underset{ KgH\in G/H/ K}{ \sum}{\rm ind}_{c_{g}: H\cap g^{-1}Kg\to K} \circ {\rm res}_{H}^{H\cap g^{-1}Kg},$$ 
  where $c_{g}$ denotes  conjugation  with  $g$. 
 \end{itemize}
 
 \end{definition}

The  following  result  was  proved in \cite{lueckequivariantcohmological},  Theorem  5.2 in page 1046. 
\begin{theorem}[Injectivity  and Mackey  functors]
Let  $G$  be  a  Group  and let $R$  be  a  commutative  ring   such  that  the  order  of  every  finite  subgroup  is   invertible in $R$.  Assume  that   $M$   is  a  mackey  functor. 

Suppose  that  the $R[W_{G}(H)]$-module  $T_{G/H}M$ is  injective  as  a $R[W_{G}(H)]$-module  for  each object. Then, $M$  is   injective  as   a $\mathcal{S}_{G}$-module, and  the  map  $\nu$  is bijective. 
\end{theorem}

A  finer structure  occurs  for cohomology  with  coefficients  in  modules  over  the  Green  functor  of  the  rational  representation  ring.  The  following  corollary  is   even  true  for  Bredon  cohomology  with  coefficients  in  such  modules.  See  \cite{lueckcreelle}, sections  6  and  7.

The   following  theorem was  proved as  6.3 in \cite{lueckcreelle}, page  221  for  Bredon  homology

\begin{theorem}
Let  $(X,A)$ be a  proper  $G$-CW  pair. Let  $M$   be  a  Mackey  functor with  module  structure  over  the Green  ring  of  rational  representations. Then, there  exists  a  decomposition

$$ H_{p}^{\mathcal{O}_{G}}(X,A) \cong\underset{H\in I}{\bigoplus}H_{p}(X^{H},A^{H}/ C_{G}(H) )\otimes_{ R[W_{G}(H)]}S_{G/H} M,   $$
where  $I$  denotes  the  set of  conjugacy  classes  of  finite  subgroups.       
      
 \end{theorem}

We  quote  now  the  most  complete  result  which  uses  the Module  structure  over  the  Green  Ring of  the  rational  representation  ring. This  appeared  as  Theorem  0.2 in  \cite{lueckcreelle}. 
 
\begin{theorem}
Let  $M$  be  a  Mackey  functor  which  admits  a  module  structure  over  the  Green  functor  of  rational  representations. For  any  group  and  any  $G$-CW  pair $(X,A)$ there  exists  a  direct  sum  decomposition 
$$H_{p}^{\mathcal{O}_{G}}(X,A) \cong\underset{H\in I}{\bigoplus}H_{p}((X^{C},A^{C})/ C_{G}(H))\otimes R[W_{G}(C)]\theta_{C}^{C}M ). $$   

Here $\theta_{C}^{C}M$ denotes  mutiplication  with  an  idempotent  in  the rational representation  ring, and  this  image  equals: 
$$ \coker \underset{D\subsetneq C}{\bigoplus} M(D)\overset{{\rm ind}_{D}^{C}}{\to} M(C).  $$

\end{theorem}  

The  explicit  use  of  such  idempotents  plays  a  role   in  delocalized  Chern  characters \cite{matthey}.

\subsection{Computations Based  on  Elementary Homological  Algebra over the  Orbit  Category}

The  fact  that  Bredon cohomology  can  be   defined  as  a  Hom construction (limit)  obtaining a  cochain  complex from  which   Bredon  cohomology  is  obtained  as usual homology allows that usual constructions (based  on  the  existence  of  resolutions and  the  concept  of  derived  functor)  in   homological  algebra  often  have  a  generalization to  Bredon  versions.

 We  present  two instances  of  these constructions:  the  Universal  Coefficient  Theorem for  Bredon cohomology of \cite{barcenasvelasquez} and  the K\"unneth  theorem  of  \cite{sanchezcoxeter}. 
    
The  following  result  appeared  in  \cite{sanchezcoxeter}, Theorem 3.1 in page 776. The  main  hypothesis  asks for the  property  that  the evaluations of  the  Bredon  Module  are  free  modules  over a  commutative  ring. Notice  that  this  does not  mean that  the   funtor  is  free  in the  sense  defined  above.  

\begin{theorem}[K\"unneth Theorem  for  Bredon cohomology]
Let  $X$ be a $G$-CW  complex and  let $Y$ be an  $H$-CW  complex. Let  $\mathfrak{F}$  and  $\mathfrak{F}^{'}$  be    families of  subgroups of  $G$  and $Y$ containing  the  isotropy  groups  of  cells  in  $X$,  respectively $Y$. 
Assume  that  $M$  and  $N$ are covariant  Bredon  functors  defined  on  the  orbit  categories $\mathcal{O}_{G}$,  respectively  $\mathcal{O}_{H}$, with  the  property  that $M(G/G^{'})$, respectively $N(H/H^{'})$  are  free  modules   for  each pair of  objects  $G/G^{'}$, $H/H^{'}$. 
 Denote  by $\mathfrak{F}\times \mathfrak{F}^{'}$ the  family of  subgroups  of  the product  which is  given as  products  of  subgroups  of $G$ and $H$, and  let  $M\otimes N$ be  the  Bredon  module  defined  in this  category.  

Then  the  product $X\times Y$   with  the  diagonal  action is a $G\times H$-CW  complex,  and  there  exists   a short    exact  sequence 

\begin{multline*}
0\to \underset{i+j=n}{\bigoplus}H_{i}^{\mathcal{O}_{\mathfrak{F}}}(X)\otimes H_{j}^{\mathcal{O}_{\mathfrak{F}^{'}}}(Y) \\ \to H_{n}^{\mathfrak{F}\times \mathfrak{F}^{'}}(X\times Y, M\otimes N) \to \underset{i+j=n}{\bigoplus} {\rm Tor}(H_{i}^{\mathcal{O}_{\mathfrak{F}}}(X), H_{j}^{\mathcal{O}_{\mathfrak{F}^{'}}}(Y))\to 0.
\end{multline*}

\end{theorem}

For  the  Universal  Coefficient Theorem  for  Bredon  cohomology,  in addition  to the  freenes  of  the  evaluation  of  the  functor on each object,  there  is  a  requirement  of  a  basis  compatibility  in  a  dual  basis  which  is a direct consequence  of  Frobenius reciprocity  for  the   complex  representation  ring  with the  characters  as  basis. We  give  the  definition  below. 

\begin{condition}\label{conditionD}
Let $G$ be  a  discrete  group, Let $ M_?$ and  $M^?$ be covariant,  res\-pec\-ti\-vely  contravariant  functors  defined  on a  subcategory $\OO$  of  the  orbit  category $\mathcal{O}$ agreeing  on  objects.  Suppose  that 

\begin{itemize}
\item There  exists  for  every  object $G/H$ a  choice of  a  finite  basis  $\{\beta_{i^ {H}}\}$    expressing $M_?(G/H)= M^ ?(G/H) $ as  the   finitely  generated,  free abelian  group on $\{\beta_{i^ {H}}\}$   and  isomorphisms  $a_{H}: M^{?}(G/H) \overset{\cong}{\rightarrow}\mathbb{Z} [ \{\beta_{i^ {H}}\}] \overset{\cong}{\leftarrow}M_{?}(G/H) :b_{H}$.

\item  For  the  covariant  functor  $\widehat{M}:=Hom_ {\mathbb{Z}}(M^ ?(\quad), \mathbb{Z})$,    the  dual  basis  $\{\widehat{\beta}_{i^ {H}}\}$ of  $Hom_\mathbb{Z}(\mathbb{Z} [ \{\beta_{i^ {H}}\}] , \mathbb{Z})$  and   the   isomorphisms $a_{H}$ and  $b_{H}$, the  following  diagram  is commutative: 
$$\xymatrix@d{ \widehat{M}(G/H) \ar[r]^{\widehat{M}(\phi)}  &  \widehat{M}(G/K) \\ \mathbb{Z}[\{\widehat{\beta}_{i^ {H}}\}] \ar[u]^{\widehat{a_{H}}} &   \mathbb{Z}[\{\widehat{\beta}_{j^ {K}}\}] \ar[u]_{\widehat{a_{K}}}  \\ \mathbb{Z} [ \{\beta_{i^ {H}}\}] \ar[u]^ {D_H}&  \mathbb{Z} [ \{\beta_{j^ {K}}\}] \ar[u]_{D_K}\\        M_{?}(G/H) \ar[r]_{M_?(\phi)} \ar[u]^{b_{H}}  & M_{?}(G/K) \ar[u]_{b_K} }$$

Where  $D_H$, $D_K$  are  the  duality  isomorphisms  associated  to  the  bases and   $\phi:G/H\to  G/K$ is  a  morphism  in the  orbit  category.

\end{itemize}
 
\end{condition}
Conditions  \ref{conditionD}  are  satisfied  in some  cases:

\begin{itemize}
\item Constant  coefficients $\mathbb{Z}$. 
\item  The  complex representation  ring  functors  defined  on the  family $\mathcal{FIN}$  of  finite  subgroups,  $\mathcal{R}^?$, $\mathcal{R}_?$. A  computation using  characters as  bases  and  Frobenius  reciprocity yields  conditions \ref{conditionD}.
\item Consider  a   discrete group $G$ and  a    normalized torsion cocycle 
 $$\alpha \in Z^2 (G,S^1),$$ take the  $\alpha$- and  $\alpha^{-1}$ twisted  representation  ring functors $ \mathcal{R}^{\alpha}_?$  $\mathcal{R}{_{\alpha}}^ ?$ defined  on the  objects  $G/H$, where  $H$   belongs  to  the family $\calfin$ of  finite  subgroups. Consider  for  every  object  $G/H$  the  cocycles $i_{H}^{*}(\alpha)$, where  $i_{H}:H\to G$ is  the  inclusion, and  assume  without loss  of  generality  that  they  are  normalized and  correspond to a  family  of  Schur  covering  groups in central extensions   $1\to \mathbb{Z}/n_{H} \to H^* \to  H\to  1 $.  

We  select  the set  $\{\beta_{H}\}$  given   as the  set  of  characters of irreducible  representations  of   $H^*$  where  $\mathbb{Z}/n_{H}$  acts by  multiplication  with a primitive  $n_{H}$-th root  of  unity.  Given  a  choice   of  sections   for  the quotient maps  $H^ *\to H$, one  can  construct  isomorphisms $^{i^{*}(\alpha)} {\mathcal{R}}      (G/H) \overset{\cong}{\rightarrow} \mathbb{Z} [ \{\beta_{H}\}] $. The  orthogonality  relations and Frobenius  reciprocity  for their  twisted characters guarantee that conditions \ref{conditionD} yield.

\end{itemize}

\begin{theorem}[Universal  Coefficient  Theorem for  Bredon Cohomology]\label{theoremUCT}

Let  $X$  be  a  proper,  finite    $G$-CW complex.  Let  $M^ ?$ and  $M_?$  be  a  pair  of functors  satisfying  conditions \ref{conditionD}.  Then, there  exists   a short  exact  sequence  of  abelian  groups

 $$ 0\to {\rm Ext}_{\mathbb{Z  }}  (H_{n-1}^{\mathcal{O}_G} (X, M_?), \mathbb{Z})\to  H^ {n}_{\mathcal{O}_{G}}(X,  M^?) \to {\rm Hom}_{\mathbb{Z  }}  (H_{n}^{\mathcal{O}_{G}} (X, M_?), \mathbb{Z}) \to  0 $$ 

\end{theorem}

\subsection{Computations  Based  on Structural Properties  of Orbit  Categories of  Groups}

There  exist examples  of  conditions  on a  group which have  consequences  on  the   particular  shape  that  an   orbit  category might  take. 

We  give  for  this  the  example  of  the  family  of  finite  groups of  a  group  which  is  a  central  extension by   a  finite  cyclic  group of  a  discrete  group which is classified  by a  second degree cohomology class  with coefficients  on the  finite  cyclic  group.

Under this  condition, there  exists  a  bijective  correspondence  between finite  subgroups   of   the  central  extensions  and  inverse  images  of finite  subgroups  in  the  original  group.

Having  the  aim  of  computing  Bredon  cohomology  with  coefficients in   twisted   complex representations,   there  exist an \emph{untwisting  procedure} to  change   twisted  coefficients  in  favor  of  the  extension  group,  and  a  certain class  of  representations. Let  us  recall  the  needed   definitions.

\begin{definition}\label{rep defzklineal}
Let  $ 1\to  \mathbb{Z}/n \mathbb{Z} \to  \tilde{H} \to  H\to 1$  be  a  central  extension. Let  $k$ be a natural number with $0\leq k\leq n$. Let  $V$  be  a  complex  vector  space.  A  $k$-central  representation of  $\tilde{H}$  is  a  homomorphism  $\tilde{H} \to {\rm GL}(V)$,  where  the  generator $t \in   \mathbb{Z} /n \mathbb{Z}$ acts  by  multiplication  by $e^{2\pi i k/n}$.

\begin{definition}\label{defkcentralbredonmodule}
The $k$-central  representation  group of  $\tilde{H}$, denoted  by  $R_k(\tilde{H})$, is the  Grothendieck  group  of  isomorphism  classes  of  $k$-central  representations  of  $\tilde{H}$. 
\end{definition}
The $k$-central representation group is a contravariant coefficient system. Given a central extension of discrete groups, $ 1\to  \mathbb{Z}/n \mathbb{Z} \to  \tilde{G} \to  H\to 1$, we denote by $\mathcal{R}^?$  the functor
\begin{align*}
\mathcal{R}^?:\Or_\calfin(\tilde{G})&\rightarrow \mathbb{Z}\text{ - }\MODULES\\
\tilde{G}/\tilde{H}&\mapsto R_k(\tilde{H}). 
\end{align*}
\end{definition}

The  following  theorem  appeared  as 4.4  in  page 57 of \cite{barcenasvelasquezrefrito}. It  was  originally  stated  for  the  classifying  space  for  proper  actions, but  it  holds  for  any  proper  $G$-CW  complex.

It  is  the  main  input  for  the  untwisting  argument  for  twisted  equivariant  $K$-Theory  of   discrete  torsion twists described  below.

\begin{theorem}\label{theo:untwisting}
Let $G$ be a discrete group and  let $\alpha\in Z^2(G;S^1)$ be a cocycle taking values in $\mathbb{Z}/n\mathbb{Z}\subseteq S^1$. Consider the extension associated to $\alpha$
$$\xymatrix{1\ar[r]&\mathbb{Z}/n\mathbb{Z}\ar[r]&G_\alpha\ar[r]^{\rho}&G\ar[r]&1.}$$

Denote  by  $X$ a  $G$-CW complex  with  finite  groups  as  cell  stabilizers.

Then, the  map $\rho$  gives  an  isomorphism  of  abelian  groups  between  the   Bredon  cohomology  groups of  $X$  with  coefficients  in  the  $\alpha $-twisted  representation  group   and  the $G_\alpha$- equivariant  Bredon  cohomology  groups  of  $X$  with  coefficients  in  the  so- called $1$-central  group representation Bredon module (defined  in  \ref{defkcentralbredonmodule}). In  symbols, 

$$H^*_{\mathcal{O}_{G}}(X;\mathcal{R}_\alpha^G)\xrightarrow{\rho^*}H^*_{\mathcal{O}_{G_\alpha}}(X;\mathcal{R}_1^{G_\alpha})$$
is an isomorphism. 
\end{theorem}

Further instances  of computations  based  on   knowledge about  the  orbit  category  of specific  examples  are  often   stated  in  terms  of  the  existence of  particular  model  for  classifying  spaces. 

These  computations hold  more  generally,  for  any equivariant  cohomology  or  homology   theory, including  Bredon  cohomology,  and  we  mention the   following    instances:

\begin{itemize}
\item   Conditions  $M$ and  $NM$  of  page  294 in  \cite{lueckclassifying}, which  hold together  for  Fuchsian  groups,  One  relator  groups, and  extensions $1\to \mathbb{Z}^{n}\to G\to F $ for  finite $F$  acting freely  outside 0. 
\item The  computation  of  equivariant  homology  theories  for classifying  spaces  of  families  of Graph product  groups  of  \cite{kasprovskililueck}.
\item Condition $C$ of  \cite{antolinflores}.

\end{itemize}

\section{Bredon Cohomology  as  Recipient  for Equivariant  Chern Characters  }

An  equivariant  Chern  Character  is a natural  transformation  between  equivariant  cohomology  theories. Let  us   briefly  recall  the  notion  of  an  equivariant  cohomology  theory.  We  stress  that   the  equivariant  cohomology  theories  considered  here  are  \emph{naive}, and remit  to  \cite{degrijsehausmannlueckpatchkoriaschwede},

\begin{definition}
 Let  $G$  be  a  group  and  fix  an  associative  ring  with  unit $R$. A $G$-Cohomology  Theory  with  values  in  $R$-modules is  a  collection  of  contravariant  functors $\mathcal{H}^{n}_{G}$ indexed  by  the  integer  numbers $\mathbb{Z}$ from  the  category  of  $G$-$CW$  pairs together  with  natural  transformations $\partial^{n}_{G}: \mathcal{H}^{n}_{G}(A):=\mathcal{H}^{n}_{G}(A, \emptyset)\to \mathcal{H}^{n+1}_{G}(X,A)  $, such  that  the  following axioms  are  satisfied: 

\begin{enumerate}
 \item{If  $f_{0}$  and  $f_{1}$  are  $G$-homotopic  maps $(X,A)\to (Y,B)$ of   $G$-CW  pairs, then  $\mathcal{H}^{n}_{G}(f_{0})=\mathcal{H}^{n}_{G}(f_{1})$ for  all n. }
 \item{Given a  pair $(X,A)$  of  $G$-$CW$  complexes, there  is  a  long  exact  sequence 
\begin{multline*}
$$ \ldots \overset{\mathcal{H}^{n-1}_{G}(i)} {\rightarrow}   \mathcal{H}^{n-1}_{G}(A)   \overset{\partial^{n-1}_{G}} {\rightarrow}   \mathcal{H}^{n}_{G}(X,A) \overset{\mathcal{H}^{n}_{G}(j)}{\rightarrow} \mathcal{H}^{n}_{G}(X) \\ \overset{\mathcal{H}^{n}_{G}(i)} {\rightarrow} \mathcal{H}^{n}_{G}(A) \overset{\partial^{n}_{G}} {\rightarrow} \mathcal{H}^{n+1}_{G}(X,A)  \overset{\mathcal{H}_{n+1}(j)}{\rightarrow} \ldots $$
\end{multline*}

 where $i:A\to X$ and  $j:X\to (X,A)$ are  the  inclusions. }

\item{Let  $(X,A)$ be  a  $G$-$CW$ pair  and  $f: A\to B$  be  a   cellular  map. The  canonical  map $(F,f): (X,A)\to (X\cup_{f} B, B)$ induces  an  isomorphism 
$$ \mathcal{H}^{n}_{G}(X\cup_{f}B, B) \overset{\cong}{\to} \mathcal{H}^{n}_{G}(X,A)$$ }

\item{ Let $\{ X_{i}\mid i\in \mathcal{I} \}$  be  a  family  of  $G$-$CW$-complexes and  denote  by $j_{i}: X_{i}\to  \coprod_{i\in \mathcal{I}} X_{i}$ the  inclusion  map. Then  the  map  
$$\Pi_{i\in \mathcal{I}}\mathcal{H}^{n}_{G}(j_{i}):  \mathcal{H}^{n}_{G}(\coprod_{i}X_{i})   \overset{\cong}{\to}  \Pi_{i\in \mathcal{I}}\mathcal{H}^{n}_{G}(X_{i})$$
is  bijective  for  each  $n\in \mathbb{Z}$. }

\end{enumerate}
A $G$-Cohomology  Theory  is  said  to  have  a  multiplicative  structure  if   there  exist  natural, graded  commutative $\cup$- products  

$$\mathcal{H}^{n}_{G}(X,A)\otimes \mathcal{H}^{m}_{G}(X,A) \to \mathcal{H}^{n+m}_{G}(X,A)$$

Let $\alpha:H\to G$ be  a  group  homomorphism and  $X$ be  a $H$-CW  complex. The  induced space ${\rm ind}_{\alpha}X,$ is defined  to be  the  $G$-CW complex  defined  as  the  quotient space $G\times X $ by  the  right  $H$-action  given  by $(g,x)\cdot h =( g\alpha(h),h^{-1}x)$. 

An Equivariant  Cohomology  Theory consists  of a  family of $G$-Cohomology Theories  $\mathcal{H}^{*}_{G}$  together  with
an induction structure determined  by  graded ring  homomorphisms 

$$ \mathcal{H}^{n}_{G}({\rm ind}_{\alpha}(X,A))\to   \mathcal{H}_{H}^{n}(X,A) $$ 
which  are  isomorphisms  for  group  homomorphisms $\alpha: H\to G$  whose kernel acts freely
on $X$ satisfying  the  following  conditions: 
\begin{enumerate}
 \item{For  any  $n$, $\partial^{n}_{H}\circ {\rm ind}_{\alpha}= {\rm ind}_{\alpha}\circ \partial^{n}_{G}$.}
 \item{For any  group  homomorphism $\beta: G\to K$ such  that $\ker \beta\circ \alpha$ acts freely on  $X$,  one  has 
$${\rm  ind}_{\alpha\circ \beta}= \mathcal{H}^{n}_{K}(f_{1}\circ {\rm ind}_{\beta}\circ{\rm ind}_{\alpha}): \mathcal{H}_{K}^{n}({\rm ind}_{\beta\circ \alpha}(X,A))\to  \mathcal{H}^{n}_{H}(X,A)$$
where $f_{1}: {\rm ind}_{\beta}{\rm ind}_{\alpha}\to {\rm ind}_{\beta\circ\alpha}$ is  the  canonical $G$-homeomorphism.}                                      

\item{For  any $n\in \mathbb{Z}$, any  $g\in G$ , the  homomorphism 

$${\rm ind}_{c_(g):G\to G}: \mathcal{H}^{n}_{G}(\rm ind)_{c(g):G\to G}(X,A))  \to \mathcal{H}^{n}_{G}(X,A) $$

agrees  with  the  map  $\mathcal{H}^{n}_{G}(f_{2})$, where  $f_{2}: (X,A)\to {\rm ind}_{c(g):G\to G}$ sends $x$  to  $(1,g^{-1}x)$ and $c(g)$  is  the  conjugation  isomorphism in $G$.}
\end{enumerate}
\end{definition}

\begin{example}[Examples  of Equivariant  Cohomology Theories]
We  now  describe  the cohomology  theories  which will be  relevant  for  the  computations  below.

\begin{enumerate}
\item Complex  Equivariant  $K$-theory  was  defined  via  vector  bundles for  finite proper  $G$-CW  complexes in \cite{oliverlueck}.  For  any proper orbit  $G/H$ one  has 
$$KU^{*}(G/H)= \begin{cases} \mathcal{R}_{\mathbb{C}}(H) \text{ $*=2k$}\\ 0 \text{$*=2k+1$}\end{cases}. $$

\item Complex, twisted  equivariant $K$-Theory with  a twist given  by  a  torsion  element  in $H^{3}(BG, \mathbb{Z})$  was  defined   in \cite{dwyer}. This  is  an  equivariant  cohomology  theory which  is  a submodule over  untwisted,  complex  $K$- theory  in  the  sense  of Oliver  and  L\"uck  described above. 

Twisted  equivariant  $K$- Theory  for  any  twist  in $H^{3}(EG\times_{G}X, \mathbb{Z})$  was  defined in  the  Fredholm  picture   in \cite{barcenasespinozajoachimuribe}. The aproppriate  axiomatic  for   twisted  equivariant  $K$-theory  for  any  third  cohomology  twist  is  that  of  parametrized  cohomology  theories \cite{barcenasespinozauribevelasquez} and  are  more  general  than  the viewpoint  adopted  here. The  equivariant  Chern Character, however,  is not  a  rational  isomorphism in the  twisted  case.

\end{enumerate}

\end{example}

\subsection{The  Atiyah-Hirzebruch Spectral  Sequence}

The Atiyah-Hirzebruch  spectral  sequence for  equivariant  cohomology  theories  was developed  by  Davis  and  L\"uck  in  \cite{davislueck}. A detailed  deduction  and  a  presentation  of the  relevant    details  is  available  in \cite{degrijsehausmannlueckpatchkoriaschwede},  page 108 construction  3.214.

\begin{theorem}
Let $\mathcal{H}^{*}$ be  an equivariant  cohomology  theory. Then,  there  exists  a spectral  sequence  which  has  $E_{2}$-term Bredon cohomology with  coefficients  in  the  functor 
$$\mathcal{H}^{q}:G/H\mapsto \mathcal{H}^{q}_{G}(G/H) $$

$$E_{2}^{p,q}= H^{q}_{\mathcal{O}_{G}} (X, \mathcal{H}^{q}), $$
which  converges contidionally  to  the equivariant  cohomology theory  modules 
$$\mathcal{H}^{*}_{G}(X). $$
\end{theorem}

It  is  a  consequence   of  the  existence  of  the equivariant Chern  character,  that  the  Atiyah-Hirzebruch  spectral  sequence rationally collapses.

\begin{remark}[The $RO(G)$-graded  Atiyah-Hirzebruch  Spectral Sequence]

There  is  a  discussion  of a  spectral  sequence  to  compute   $RO(G)$- graded cohomology  theories  out  of  $RO(G)$-graded Bredon  cohomology  in \cite{kronholm}.  
\end{remark} 

\begin{remark}[The p-Chain Spectral  Sequence]
An  alternative  to  the  Atiyah-Hirzebruch  spectral  sequence  is  given  by  the  $p$- Chain  spectral  sequence  of \cite{davislueckp}. The  spectral  sequence  is  constructed  in  the  setting  of  spaces  over  a  category \cite{davislueck}  and   it  is   thus  valid  for  any  equivariant  cohomology  theory, particularly  Bredon cohomology. 
\end{remark}

\begin{remark}[Third Differential of  the Equivariant  Atiyah-Hirzebruch  Spectral Sequence   for  Equivariant  Complex  $K$-Theory]

Let  us   restrict  to  complex  equivariant  $K$- theory.  While  the equivariant  Atiyah-Hirzebruch  spectral  sequence  rationally  collapses, and   for  non-equivariant complex  $K$-theory there  exist closed formulas  for  the  first  non-vanishing  differential, $d_{3}$ (and  even  for  higher degree  ones,  in terms  of  secondary  cohomology  operations),  as  of  2022  there  exists  no  closed  formula  for  the  third  differential 
$$d_{3}: H^{p}_{\mathcal{O}_{G}}(X, \mathcal{R}_{\mathbb{C}})\to H^{p+3}_{\mathcal{O}_{G}}(X,\mathcal{R}_{\mathbb{C}}). $$

See  \cite{degrijseleary} for  a  discussion  of  the failure  of  the  integral  differential  to be  an  isomorphism. 

Particular  instances  of  the spectral sequence (without exhausting  the  most  general  closed  formula) have  been  discused  in  the  literature. 

Work  by  Uribe  and G\'omez \cite{uribegomez} introduced  a decomposition  of  (non twisted)  equivariant complex  $K$- theory of  finite  groups  $G$ which  have  a normal  abelian subgroup $A$ which  acts  trivially  on   a  finite  $G$-CW  complex $X$ in  summands  of  twisted  equivariant  $K$- Theory  corresponding  to  a  number  of  twists corresponding  to  irreducible  representations  of $A$.

The  outcome  is  that, under  these  additional  conditions, which do  not  hold for  a  general action  of  a finite group $G$ on  a $G$-CW  complex $X$,   they  are  able  to  identify the  third  differential  of  the (untwisted) equivariant Atiyah-Hirzebruch  spectral  sequence with a  special  instance  of  the  third  differential  of  the  twisted Segal spectral  sequence  constructed  in \cite{barcenasespinozauribevelasquez}.  

While the  differentials  of  the equivariant  Atiyah-Hirzebruch spectral  sequence  are    natural  transformations,  and  even  in the  parametrized  setting  they are  identified  by Theorem 5.7 in  \cite{barcenasbrown} in homotopy theoretical  terms as  maps  between  classifying  spectra  for  Bredon  cohomology, and  there  exist constructions of  the  cohomology  operations  in \cite{ginot}, it  has   not  been  possible  to  give a  complete  list  of  candidates  for  the  relevant  cohomological  operations between  Bredon  cohomology  groups. 

In the  $RO(G)$-graded  setting, the  definition  of  the adequate  version  of  the  Steenrod  Algebra  goes  back  to Oru\c{c} \cite{oruc}, and  efforts  to address  the  analogous  problem  of   determining  the  possible  operations   are \cite{ricka},  and  \cite{sankar}.

\end{remark}

\subsection{Equivariant  Cohomological Chern Characters}

The  equivariant  Chern  character for was addressed  first  by Slominska  for  equivariant $K$-Theory  of  finite  groups in \cite{slominska76}.

In the  context  of the Baum-Connes  conjecture, the  need  for  decomposition of equivariant  $K$-homology into  informations  of  fixed  point  sets  of  finite  cyclic  groups led  to a  specific  construction,  named  the  delocalized  Chern  character  in \cite{baumconneschern},  and  the  formalization  in  the  terms  refered  here  was  done  by  L\"uck   mainly  in the  articles \cite{lueckcreelle}, \cite{lueckequivariantcohmological}.

The  following  result  was  proved  as  Theorem  4.2  in  page 1041 of  \cite{lueckequivariantcohmological}. 

\begin{theorem}[The  Equivariant  Chern Character]
Let  $R$ be  a  ring  containing  the  rational  numbers.  Let $\mathcal{H}^*$ be  a  proper  equivariant  cohomology  theory  with values in $R$-modules.  Suppose  that  the $\mathcal{S}_{G}$-module $\mathcal{H}^q\circ {\rm pr}$ is injective  as  $\mathcal{S}_{G}$-module for  every  group $G$ and  every $q\in \mathbb{Z}$.  Then, we obtain a transformation  of  proper  equivariant  cohomology  theories 
$$ {\rm ch}^{n}: \mathcal{H}^{n}_{?}\longrightarrow \underset{p+q=n}{\prod} H ^{p}_{\mathcal{S}_{G}}(X, \mathcal{H}^{q}).$$ 
 The  $R$- map  is bijective  for  all proper  relatively  finite  $G$-CW  pairs $(X,A)$.  if  $\mathcal{H}^*$ satisfies  the  disjoint  union  axiom,  then  the $R$- map  is  bijective  for  all proper  $G$-CW  pairs  $(X,A)$. 
 \end{theorem}
 The  natural  transformation  is  constructed  using  a composition  of  eight different  construction in  page  1040  of  L\"uck,  to  where  we  refer  for  further  details.

\begin{remark}[Orbifold  Version]

Adem and  Ruan introduced  an  orbifold  bersion of   both   twisted  complex  $K$-theory and  Bredon  cohomology for  a discrete torsion  twist.

 The latter one  turns out  to  be  isomorphic  to Chen-Ruan  cohomology  of  orbifolds in \cite{ademruan}.  See \cite{ademruanleida},   3.3  and 3.10  in  pages 60,  and 77  for  a more  detailed  exposition.
  
\end{remark}

\begin{remark}[Delocalized  Version]
In  connection  with  the rationalized Baum-Connes  assembly  map, delocalization refers  to modified  versions  of  the  equivariant Chern  character, which  within a  geometric  setting  can  be  thought  of  being  defined before inverting  the  Thom  class  of  normal  bundles  of  inclusions  of  fixed  point  sets, \cite{baumconneschern}, paragraph 9. A  good  discussion of  these versions  of  the  Chern Character, including  the   relation  to  the equivariant  Chern  character  presented  here  is  given  in \cite{matthey}.   

\end{remark}

\begin{remark}[Homotopy Theoretical refinements]
As  of 2022,  the  most  refined  homotopy theoretical  versions  of  the  Chern  Character  are   presented  in   \cite{lackmann}.  The  results  initiated  in the  Author's  Ph. D Thesis  and  include  a  study  of   external  duality,  and  a  homology  representation theorem.  See  also  \cite{lackmannli}. 
\end{remark}

\bibliographystyle{amsplain}
\bibliography{bredon}

\end{document}